\documentclass[12pt,a4paper]{article}

\usepackage{graphicx}
\usepackage{amsmath,amsfonts}
\usepackage{mathrsfs}
\usepackage{authblk}
\newtheorem{theorem}{Theorem}

\newtheorem{proposition}{Proposition}

\let\scr\mathscr

\def\Pb{\mathbf{P}}

\def\Ex{\mathbf{E}}

\def\BB{\mathbb{B}}
\def\KK{\mathbb{K}}
\def\EE{\mathbb{E}}

\def\1{\mbox{1\hspace{-.25em}I}}
\def \Liminf{\mathop{\underline{\lim}}\limits}

\begin{document}
\title{On Parameter Estimation  of the Hidden Gaussian Process in perturbed
 SDE.}
\author[1]{Yu. A. Kutoyants}
\author[2]{L. Zhou}
\affil[1]{Le Mans University, Le Mans, France}
\affil[1]{Tomsk State University, Tomsk, Russia}
\affil[2]{School of Mathematics and Statistics, Shandong University, Weihai, China}

\date{}

\maketitle
\begin{abstract}
We present results on parameter estimation and non-parameter estimation
of the linear partially observed
Gaussian system of stochastic differential equations. We propose new one-step
estimators which have the same asymptotic properties as the MLE, but much more
simple to calculate, the estimators are
so-called "estimator-processes".  The construction of the estimators is based
on the equations of Kalman-Bucy filtration and the asymptotic corresponds to
the {\it small noises} in the observations and state (hidden process)
equations. We propose conditions which provide the consistency and asymptotic
normality and asymptotic efficiency of the estimators.

\end{abstract}


\section{Introduction}

Let us consider the problem of parameter estimation for partially observed
linear system. The observed process is:
\begin{align}
\label{1}
{\rm d}X_t=f(\theta,t)Y_t\,{\rm d}t+\varepsilon\sigma\left(t\right){\rm
  d}W_t,\quad X_0=0,
\end{align}
where the {\it  hidden} process $Y_t$ is  solution  of the equation
\begin{align}
\label{2}
{\rm d}Y_t=a(\vartheta,t)Y_t\,{\rm d}t+\varepsilon b\left(t \right){\rm d}
V_t,\quad Y_0=y_0\neq 0.
\end{align}
Here $W_t, 0\leq t\leq T$ and $V_t, 0\leq t\leq T$ are independent Wiener
processes and the functions $f\left(\cdot \right), \sigma \left(\cdot
\right),a\left(\cdot \right), b\left(\cdot \right) $ are known. The unknown
parameter $\vartheta \in \Theta =\left(\alpha ,\beta \right)$, $\left|\alpha
\right|+\left|\beta \right|< \infty $.

Therefore we have to estimate $\vartheta $ by observations
$X^T=\left(X_t,0\leq t\leq T\right)$. The properties of the estimators we
describe in the asymptotic of small noises ($\varepsilon \rightarrow 0$) in
the observation \eqref{1} and state \eqref{2} equations. The maximum
likelihood estimator (MLE) $\hat\vartheta _\varepsilon $ and Bayes estimator
(BE) $\tilde\vartheta _\varepsilon $ under regularity conditions are
consistent, asymptotically normal
\begin{align}
\label{3}
\frac{\hat\vartheta _\varepsilon-\vartheta _0}{\varepsilon }\Longrightarrow
     {\cal N}\left(0, {\rm I}\left(\vartheta _0\right)^{-1}\right),\qquad
     \frac{\tilde\vartheta _\varepsilon-\vartheta _0}{\varepsilon
     }\Longrightarrow {\cal N}\left(0, {\rm I}\left(\vartheta
     _0\right)^{-1}\right)
\end{align}
and asymptotically efficient \cite{Kut84}. Here $\vartheta _0$ is the true
value and ${\rm I}\left(\vartheta _0\right) $ is Fisher information. The
construction of these estimators is based on the likelihood ratio function
\begin{align*}
L\left(\vartheta ,X^T\right)=\exp \left\{\int_{0}^{T}\frac{f\left(\vartheta
  ,t\right)m\left(\vartheta ,t\right)}{\varepsilon ^2\sigma
  \left(t\right)^2}{\rm d}X_t -\int_{0}^{T}\frac{f\left(\vartheta
  ,t\right)^2m\left(\vartheta ,t\right)^2}{2\varepsilon ^2\sigma
  \left(t\right)^2}{\rm d}t\right\}
\end{align*}
and is given by the following relations
\begin{align}
\label{4}
L\left(\hat\vartheta _\varepsilon,X^T\right)=\sup_{\vartheta \in\Theta
}L\left(\vartheta ,X^T\right),\quad \qquad \tilde\vartheta _\varepsilon
=\frac{\int_{\Theta }^{}\vartheta p\left(\vartheta \right)L\left(\vartheta
  ,X^T\right){\rm d}\vartheta }{\int_{\Theta }^{}p\left(\vartheta
  \right)L\left(\vartheta ,X^T\right){\rm d}\vartheta}.
\end{align}
Here $m\left(\vartheta ,t\right)=\Ex_\vartheta \left(Y_t|X_s,0\leq s\leq
t\right) $ is conditional expectation, satisfying the Kalman-Bucy filtration
equations
\begin{align}
\label{5}
{\rm d}m\left(\vartheta ,t\right)&=a\left(\vartheta
  ,t\right)m\left(\vartheta
,t\right){\rm d}t\nonumber\\
&\qquad \qquad \qquad +\frac{\gamma \left(\vartheta ,t\right)f\left(\vartheta
  ,t\right) }{\varepsilon ^2\sigma \left(t\right)^2} \left[{\rm
  d}X_t-f\left(\vartheta ,t\right)m\left(\vartheta ,t\right){\rm d}t \right],\\
 \frac{\partial \gamma \left(\vartheta ,t\right)}{\partial
  t}&=2a\left(\vartheta ,t\right)\gamma \left(\vartheta ,t\right)-\frac{\gamma
  \left(\vartheta ,t\right)^2f\left(\vartheta ,t\right)^2}{\varepsilon
  ^2\sigma \left(t\right)^2}+\varepsilon ^2b\left(t\right)^2,
\label{6}
\end{align}
with initial values $m\left(\vartheta ,0\right)=y_0 $ and $\gamma
\left(\vartheta ,0\right)=0 $.  Recall that $ \gamma \left(\vartheta
,t\right)=\Ex_\vartheta \left(m\left(\vartheta ,t\right)-Y_t \right)^2 $.
It is evident that the construction of the MLE and BE according to  the relations
 \eqref{4} is computationally hard problem because we need solutions of the
 system \eqref{5}-\eqref{6}   {\it for all} $\vartheta \in\Theta $. Therefore
 numerical realization of these construction is difficult to do.

The problem of construction of adaptive Kalman filter was treated by many
authors in engineering literature, where the models of observations are mainly
of discrete time form (see, e.g., \cite{Gus00}, \cite{HMS03},\cite{HCCL03}, \cite{Ru91} and
references therein). The identification of continuous time partially observed
systems were studied as well by some authors (see, e.g., \cite{A83},
\cite{CMT05}, \cite{EAM95}, \cite{EM02}, \cite{Kut84}, \cite{Kut04},
\cite{K19}, \cite{Kut19} and
references therein).  Similar problems for continuous time hidden telegraph
process were studied in \cite{PCh09}, \cite{KhK18}.

In the present work our goal is to propose other estimators, which have the
same asymptotic properties as mentioned in \eqref{3} but these estimators can
be much more easily calculated. Moreover, we construct {\it estimator-process},
i.e., estimator which evaluate with time. These constructions are based on
the One-step score-function approach, which consists in two-step estimation
procedure. This means, firstly in using a small part of initial observations,
we obtain a consistent estimator  of the unknown
parameter, and then with the help of this estimator and score-function we
construct the One-step MLE and One-step MLE-process.  This approach in general
is well
known. First such one-step procedure was proposed by Fisher \cite{F25}. Then
it was used  by many authors, see, e.g., \cite{G84},
\cite{KhK18},\cite{KU15},\cite{KM16},\cite{Kut17}.

  The construction is done in two steps. First we obtain a
preliminary estimator $\bar\vartheta_{\tau_\varepsilon}$ by the observation
$X^{\tau_\varepsilon}=(X_t,\, 0\leq t\leq \tau_\varepsilon)$, where
$\tau_\varepsilon=\varepsilon^\delta\rightarrow0$. Then using this estimator
and one-step MLE structure we obtain the estimator
$$
\vartheta_{t,\varepsilon}^*=\bar\vartheta_{\tau_\varepsilon}+
I_{\tau_\varepsilon}(\bar\vartheta_{\tau_\varepsilon})^{-1}
\int_{\tau_\varepsilon}^t\frac{\dot
  M(\bar\vartheta_{\tau_\varepsilon},s)}{\sigma\left(s\right)^2}\left[{\rm d}
  X_s-f(\bar\vartheta_{\tau_\varepsilon},s)m(\bar\vartheta_{\tau_\varepsilon},s){\rm d}
  s\right]
$$ where dot means the derivation w.r.t. $\vartheta$, $m(t, \vartheta)$ is
the conditional expectation of $Y_t$ w.r.t. $\{X_s,\,0\leq s\leq t\}$ and
$M\left(\vartheta ,s\right)= f\left(\vartheta ,s\right)y\left(\vartheta
,s\right)$.

\section{Preliminary estimator}

We say that $h\left(\vartheta ,t \right)\in {\cal C}_b$ if the  function
$h\left(\vartheta ,t \right), \vartheta \in\Theta , t\in \left[0,T\right]$ is bounded;
$h\left(\vartheta ,t \right)\in {\cal C}^{\left(1\right)}_\vartheta $ or
$h\left(\vartheta ,t \right)\in {\cal C}^{\left(1\right)}_t$ if the function
$h\left(\cdot ,\cdot \right)$ is continuously differentiable on $\vartheta $
or $t$ respectively. The derivatives w.r.t. $\vartheta $ we denote $\dot
h\left(\vartheta ,t\right)$ and the derivatives w.r.t. $t$ we write as
$h'\left(\vartheta ,t\right) $.

Let us introduce the {\it Conditions} ${\scr R}$:

${\scr R}_1.$ {\it The functions $f\left(\cdot,\cdot  \right), \sigma \left(\cdot
\right),a\left(\cdot,\cdot  \right)  , b\left(\cdot \right)\in {\cal C}_b $.}

${\scr R}_2.$ {\it The functions $f\left(\vartheta ,\cdot  \right),
  a\left(\vartheta ,\cdot  \right) \in {\cal C}^{\left(1\right)}_\vartheta $. The function $\dot
  f\left(\vartheta ,t  \right)\in {\cal C}^{\left(1\right)}_t$.}

${\scr R}_3.$ {\it The function $f\left(\vartheta ,0  \right)$ and its
  derivative  $\dot f\left(\vartheta ,0  \right)$ are separated
  from zero }
\begin{align*}
\inf_{\vartheta \in\Theta } \left|f\left(\vartheta ,0  \right)\right|>0,\qquad
\inf_{\vartheta \in\Theta }\left|\dot f\left(\vartheta ,0  \right)\right|>0.
\end{align*}

${\scr R}_4.$ {\it The initial value $y_0\not=0$. }

Recall that  ${\scr R}_4$ is a necessary condition for existence of consistent
estimator for the model of observations \eqref{1}-\eqref{2}. If $y_0=0$, then
we can denote $\hat X_t=\varepsilon^{-1} X_t$ and $\hat Y_t=\varepsilon^{-1} Y_t$ and
rewrite the system \eqref{1}-\eqref{2} as follows
\begin{align*}
{\rm d}\hat X_t&=f(\vartheta,t)\hat Y_t\,{\rm d}t+ \sigma\left(t\right){\rm d}W_t,\quad \hat X_0=0,\\
{\rm d}\hat Y_t&=a(\vartheta,t)\hat Y_t\,{\rm d}t+ b\left(t \right){\rm d}
V_t,\quad \hat Y_0=0.
\end{align*}
Hence this system does not depend on $\varepsilon $ and the consistent
estimation is impossible. This was noted by Khasminskii \cite{Kh05}. The condition
which makes main sense here is ${\scr R}_3$. Below we consider how can be
constructed preliminary estimator if this condition is replaced by another
one.  Without
loss of generality we suppose that $y_0>0, f\left(\vartheta ,0  \right)>0 $  and
$\dot f\left(\vartheta ,0  \right)>0 $.

Let us denote $x_t(\vartheta)$ and $y_t(\vartheta) $ the solutions of the
equations  \eqref{1}, \eqref{2}  for $\varepsilon =0$:
\begin{align*}
\frac{\partial  x_t(\vartheta)}{\partial  t}&=f(\vartheta,t)
y_t(\vartheta),\qquad x_0\left(\vartheta \right)=0, \\
\frac{\partial  y_t(\vartheta)}{\partial
  t}&=a(\vartheta,t)y_t(\vartheta),\qquad y_0(\vartheta)= y_0.
\end{align*}
Hence
\begin{align*}
x_t(\vartheta)&=y_0\int_0^tf(\vartheta,s) A\left(\vartheta ,s\right)   {\rm d}s,\\
y_t(\vartheta)&=y_0\; A\left(\vartheta ,t\right),\quad
A\left(\vartheta ,t\right)=\exp\left\{\int_{0}^{t}a\left(\vartheta ,r\right){\rm
  d}r\right\}.
\end{align*}
The true value we denote as $\vartheta _0\in\Theta $. We have equalities
\begin{align*}
Y_t-y_t\left(\vartheta _0\right)&=\varepsilon
\int_{0}^{t}A\left(\vartheta _0,t,s\right)b\left(s\right){\rm d}V_s, \qquad
A\left(\vartheta _0,t,s\right)=\frac{A\left(\vartheta
  _0,t\right)}{A\left(\vartheta _0,s\right)} ,\\
X_t-x_t\left(\vartheta _0\right)&=\varepsilon \int_{0}^{t}f\left(\vartheta
_0,s\right)\int_{0}^{s} A\left(\vartheta _0,s,q\right)b\left(q\right){\rm d}V_q\;{\rm
  d}s+\varepsilon \int_{0}^{t}\sigma \left(s\right)\,{\rm d}W_s\\
&=\varepsilon \int_{0}^{t}b\left(q\right)\int_{q}^{t }f\left(\vartheta
_0,s\right)  A\left(\vartheta _0,s,q\right){\rm
  d}s\; {\rm d}V_q+\varepsilon \int_{0}^{t}\sigma \left(s\right)\,{\rm d}W_s.
\end{align*}
 Hence we can write
\begin{align*}
Y_t-y_t\left(\vartheta _0\right)=\varepsilon \xi _t,\qquad
X_t-x_t\left(\vartheta _0\right)=\varepsilon \eta _t ,\qquad 0\leq t\leq T,
\end{align*}
where $\xi _t$ and $\eta _t$ are Gaussian processes with $\Ex_{\vartheta_0}
\xi _t=0$, $\Ex_{\vartheta_0} \eta  _t=0$ and for any $p>0$
\begin{align}
\label{7}
\Ex_{\vartheta_0}\left|\xi _t\right|^p<C_1t^{p/2},\qquad \quad
\Ex_{\vartheta_0}\left|\eta  _t\right|^p<C_2t^{p/2}.
\end{align}
 The constants $C_1>0,C_2>0 $ do not depend on $\vartheta _0$ and
 $t\in\left[0,T\right]$ (see \cite{LS}).

  Introduce the functions
\begin{align*}
x_m\left(t\right)=\inf_{\vartheta \in\Theta }x_t \left(\vartheta
\right)<\sup_{\vartheta \in\Theta }x_t \left(\vartheta
\right)=x_M\left(t\right).
\end{align*}
In the vicinity of the point $t=0$ the function $x_t \left(\vartheta
\right)$ is monotone increasing on $\vartheta $ and we have
\begin{align*}
x_m\left(t\right)=x_t \left(\alpha \right),\qquad
x_M\left(t\right)=x_t \left(\beta \right).
\end{align*}
Further, following \cite{KhK18} and \cite{Kut19} we put $\tau _\varepsilon
=\varepsilon ^\delta,\delta >0 $,   introduce sets
\begin{align*}
&\BB_\varepsilon=\left\{x_m\left(\tau
_\varepsilon \right)<X_{\tau _\varepsilon }< x_M\left(\tau
_\varepsilon \right) \right\},\\
&\BB_\varepsilon ^-=\left\{X_{\tau _\varepsilon }\leq
  x_m\left(\tau_\varepsilon \right) \right\},\qquad \quad  \BB_\varepsilon
  ^-=\left\{X_{\tau _\varepsilon }\geq x_M\left(\tau
_\varepsilon \right) \right\} ,
\end{align*}
  and define the estimator
\begin{align}
\label{8}
\bar\vartheta _{\tau _\varepsilon }=\alpha \1_{\left\{\BB_\varepsilon
  ^-\right\}}+\mu _\varepsilon \1_{\left\{\BB_\varepsilon\right\}} +\beta
\1_{\left\{\BB_\varepsilon ^+\right\}},
\end{align}
where $\mu _\varepsilon $ is solution of the equation
$x_{\tau _\varepsilon }(\mu _\varepsilon)=X_{\tau_\varepsilon}$.

This is preliminary estimator which will be used in the next section for
construction of asymptotically efficient estimator.

\begin{theorem}
\label{T1}
Suppose that the conditions ${\scr R}$ be fulfilled and $\delta \in (0,2)$,
then uniformly on compacts $\KK\in\Theta $ the estimator $\bar\vartheta _{\tau
  _\varepsilon } $ is consistent, i.e., for any $\nu >0$ and any $\KK$
\begin{align}
\label{9}
\sup_{\vartheta _0\in \KK}\Pb_{\vartheta _0}\left(\left|\bar\vartheta _{\tau
  _\varepsilon }-\vartheta _0 \right|>\nu \right)\longrightarrow 0
\end{align}
as $\varepsilon \rightarrow 0$.  Moreover, for any $p>0$
\begin{align}
\label{10}
\sup_{\vartheta _0\in
  \KK}\Ex_{\vartheta_0}|\bar\vartheta_{\tau_\varepsilon}-\vartheta_0|^p\leq
C\left(\frac{\varepsilon}{{\sqrt{\tau_\varepsilon}}}\right)^p\longrightarrow 0.
\end{align}
\end{theorem}
{\bf Proof.} We have
\begin{align*}
\Pb_{\vartheta _0}\left(\left|\bar\vartheta _{\tau _\varepsilon }-\vartheta _0
\right|>\nu \right)&=\Pb_{\vartheta _0}\left(\left|\bar\vartheta _{\tau
  _\varepsilon }-\vartheta _0 \right|>\nu,\BB_\varepsilon ^- \right)+\Pb_{\vartheta
  _0}\left(\left|\mu  _{\varepsilon }-\vartheta _0 \right|>\nu,\BB_\varepsilon
\right)\\
&\qquad +\Pb_{\vartheta _0}\left(\left|\bar\vartheta _{\tau _\varepsilon
}-\vartheta _0 \right|>\nu,\BB_\varepsilon ^+ \right)\\
&\leq  \Pb_{\vartheta _0}\left(\BB_\varepsilon ^- \right)+\Pb_{\vartheta
  _0}\left(\left|\mu  _{\varepsilon }-\vartheta _0 \right|>\nu,\BB_\varepsilon
\right)+\Pb_{\vartheta _0}\left(\BB_\varepsilon ^+ \right).
\end{align*}

 Let $\left[\alpha _K,\beta _K\right]\subset \Theta $ be such
that $\KK\subset \left[\alpha _K,\beta _K\right] $.
We can write
\begin{align*}
\dot x_t\left(\vartheta \right)=y_0 \int_{0}^{t}\left[\dot f\left(\vartheta
  ,s\right)+ f\left(\vartheta,s\right)\int_{0}^{s}\dot a\left(\vartheta ,r\right){\rm d}r\right]
  A\left(\vartheta ,s\right){\rm d}s .
\end{align*}
By conditions ${\scr R}_2$, ${\scr R}_3$ there exists $\kappa _*>0$ such that
for sufficiently small $t$ we have the estimate
\begin{align*}
\inf_{\vartheta \in\Theta }\dot x_t\left(\vartheta \right)\geq \kappa _*t.
\end{align*}
 With the help of \eqref{7}, we have
\begin{align*}
\sup_{\vartheta _0\in \KK}\Pb_{\vartheta _0}\left(\BB_\varepsilon
  ^- \right)&=\sup_{\vartheta _0\in \KK}\Pb_{\vartheta _0}\left( X_{\tau
  _\varepsilon }-x_{\tau _\varepsilon }\left(\vartheta _0\right)\leq
  x_m\left(\tau_\varepsilon \right) -x_{\tau _\varepsilon }\left(\vartheta
  _0\right)  \right)\\
&\leq \sup_{\vartheta _0\in \KK}\Pb_{\vartheta _0}\left( \left|X_{\tau
  _\varepsilon }-x_{\tau _\varepsilon }\left(\vartheta _0\right)\right|\geq
  x_{\tau _\varepsilon }\left(\vartheta
  _0\right)-x_{\tau _\varepsilon }\left(\alpha  \right)   \right)\\
&\leq \sup_{\vartheta _0\in \KK}\Pb_{\vartheta _0}\left( \left|X_{\tau
  _\varepsilon }-x_{\tau _\varepsilon }\left(\vartheta _0\right)\right|\geq
  x_{\tau _\varepsilon }\left(\alpha _K\right)-x_{\tau _\varepsilon
  }\left(\alpha  \right)   \right)\\
&\leq \sup_{\vartheta _0\in \KK}\Pb_{\vartheta _0}\left( \left|X_{\tau
  _\varepsilon }-x_{\tau _\varepsilon }\left(\vartheta _0\right)\right|\geq
  \kappa _*\left(\alpha _K-\alpha  \right)\tau _\varepsilon \right)\\
&\leq \sup_{\vartheta _0\in \KK}\frac{\Ex_{\vartheta _0}\left|X_{\tau
  _\varepsilon }-x_{\tau _\varepsilon }\left(\vartheta
    _0\right)\right|^p}{\kappa _*^p
\left(\alpha _K-\alpha  \right)^p\tau _\varepsilon^{p}}
\leq C\;\frac{\varepsilon ^p}{\tau _\varepsilon^{p/2} }.
\end{align*}
The similar estimate we have for the probability $\Pb_{\vartheta
  _0}\left(\BB_\varepsilon ^+ \right)$. Further
\begin{align*}
\Pb_{\vartheta _0}\left(\left|\mu _\varepsilon -\vartheta _0\right|>\nu ,
\BB_\varepsilon \right)&=\Pb_{\vartheta _0}\left(\frac{\left|X_{\tau
    _\varepsilon }-x_{\tau _\varepsilon }\left(\vartheta _0\right)\right|}{ \dot
  x_{\tau _\varepsilon }\left(\tilde \mu _\varepsilon \right)}>\nu ,
\BB_\varepsilon \right)\leq C\frac{\varepsilon ^p}{\nu ^p\tau _\varepsilon^{p/2} }.
\end{align*}

For the moments  we have
\begin{align*}
\Ex_{\vartheta _0}\left| \bar\vartheta _{\tau _\varepsilon }-\vartheta
_0\right|^p&\leq \left(\vartheta _0 -\alpha \right)^p \Pb_{\vartheta
  _0}\left(\BB_\varepsilon ^- \right)+\left(\beta -\vartheta _0 \right)^p \Pb_{\vartheta
  _0}\left(\BB_\varepsilon ^+ \right)\\
&\quad +\Ex_{\vartheta _0}\left| \mu  _{\tau _\varepsilon }-\vartheta
_0\right|^p\1_{\left\{\BB_\varepsilon \right\}}
\leq C\frac{\varepsilon ^p}{\tau _\varepsilon^{p/2} }.
\end{align*}

\bigskip
 {\bf Case} $f\left(\vartheta ,t\right)=f\left(t\right)$.

 The
condition ${\scr R}_3$ is not fulfilled. Introduce another conditions

${\scr R}_2'$ {\it The function $a\left(\vartheta ,t\right)\in{\cal
  C}^{\left(1\right)}_\vartheta $  and $\dot a\left(\vartheta ,t\right)\in{\cal
  C}^{\left(1\right)}_t $ }

${\scr R}_3'$ {\it The function $a\left(\vartheta ,t\right)$  and its
  derivative $\dot a\left(\vartheta ,t\right)$  are separated from zero}
\begin{align*}
\inf_{\vartheta \in\Theta }a\left(\vartheta ,0\right)>0,\qquad \quad
\inf_{\vartheta \in\Theta }\dot a\left(\vartheta ,0\right)>0.
\end{align*}

\begin{proposition}
\label{P1}
Let the conditions ${\scr R}_1$, ${\scr R}_2'$, ${\scr R}_3'$,${\scr R}_4$ be
fulfilled and $\tau _\varepsilon =\varepsilon ^\delta $ with  $\delta \in
\left(0,2/3\right)$. Then the estimator
$\bar\vartheta _\varepsilon $ is uniformly on
compacts $\KK\subset \Theta $ consistent and
\begin{align}
\label{}
\Ex_{\vartheta _0}\left| \bar\vartheta _{\tau _\varepsilon }-\vartheta
_0\right|^2\leq C\varepsilon ^{2-3\delta }.
\end{align}
\end{proposition}
{\bf Proof.} The proof of Theorem \ref{T1} can be applied here with the only
difference in the estimation of $\dot x_t\left(\vartheta \right)$. We have
\begin{align*}
\dot x_t\left(\vartheta \right)=y_0 \int_{0}^{t} \int_{s}^{t}\dot
a\left(\vartheta ,r\right){\rm d}r
  A\left(\vartheta ,t,s\right){\rm d}s\geq \kappa ^*t^2.
\end{align*}
Hence from the proof of Theorem  \ref{T1} we obtain the estimate
\begin{align*}
\Ex_{\vartheta _0}\left| \bar\vartheta _{\tau _\varepsilon }-\vartheta
_0\right|^2\leq C\frac{\varepsilon ^2}{\tau _\varepsilon ^3}.
\end{align*}

\section{One-step MLE}

Let us re-write the equation \eqref{5}-\eqref{6} as follows
\begin{align*}
{\rm d}m\left(\vartheta ,t\right)&=\left[a\left(\vartheta
  ,t\right)-D\left(\vartheta ,t\right) f\left(\vartheta
  ,t\right)\right]m\left(\vartheta ,t\right){\rm d}t+D\left(\vartheta
,t\right) {\rm d}X_t,\\ \frac{\partial \gamma_* \left(\vartheta
  ,t\right)}{\partial t}&=2a\left(\vartheta ,t\right)\gamma_* \left(\vartheta
,t\right)-\frac{\gamma_* \left(\vartheta ,t\right)^2f\left(\vartheta
  ,t\right)^2}{\sigma \left(t\right)^2}+b\left(t\right)^2,\quad \gamma_*
\left(\vartheta ,0\right)=0,
\end{align*}
where we denoted $\gamma_* \left(\vartheta ,t\right)=\gamma \left(\vartheta
,t\right)/\varepsilon ^{2}$ and $D\left(\vartheta ,t\right)= \gamma_*
\left(\vartheta ,t\right)f\left(\vartheta ,t\right)/\sigma
\left(t\right)^{2}$.

Introduce the  conditions

$\tilde{\scr R}_2.$ {\it The functions $f\left(\vartheta ,\cdot  \right),
  a\left(\vartheta ,\cdot  \right) \in {\cal C}^{\left(2\right)}_\vartheta $. The function $\dot
  f\left(\cdot  ,t  \right)\in {\cal C}^{\left(1\right)}_t$.}

 It can be shown that the random process $m\left(\vartheta
,t\right),0\leq t\leq T $ is continuously differentiable w.r.t. $\vartheta $
with probability 1 (see, e.g.,  \cite{Kut94}). The derivative $\dot m\left(\vartheta
,t\right) $ satisfies the equation
\begin{align}
{\rm d}\dot  m\left(\vartheta ,t\right)&=\left[a\left(\vartheta
  ,t\right)-D\left(\vartheta ,t\right) f\left(\vartheta
  ,t\right)\right]\dot m\left(\vartheta
,t\right){\rm d}t+\dot D\left(\vartheta ,t\right) {\rm   d}X_t \nonumber   \\
&\quad + \left[\dot a\left(\vartheta
  ,t\right)-\dot D\left(\vartheta ,t\right)f\left(\vartheta
  ,t\right)-D\left(\vartheta ,t\right)\dot f\left(\vartheta
  ,t\right) \right] m\left(\vartheta
,t\right){\rm d}t .
\label{12}
\end{align}
Therefore for $\vartheta =\vartheta _0$ and $\varepsilon =0$ we obtain
deterministic
function $\dot y\left(\vartheta_0 ,t\right)\equiv \left.\dot m\left(\vartheta_0
,t\right)\right|_{\varepsilon =0} $, which can be written as follows
\begin{align}
\dot y\left(\vartheta_0 ,t\right)&=\int_{0}^{t}e^{B\left(\vartheta_0
  ,t,s\right)} \left[\dot a\left(\vartheta_0 ,s\right)-D\left(\vartheta_0
  ,s\right)\dot f\left(\vartheta_0 ,s\right) \right]y_s\left(\vartheta _0\right){\rm d}s,
\label{13}
\end{align}
where
\begin{align*}
B\left(\vartheta_0 ,t,s\right)=\int_{s}^{t}\left[a\left(\vartheta
  _0,r\right)-D\left(\vartheta_0 ,r\right) f\left(\vartheta_0 ,r\right)
  \right]{\rm d}r.
\end{align*}
Note that $\dot y\left(\vartheta_0 ,t\right)\not= \dot y_t\left(\vartheta
_0\right) $.

The Fisher information we define as follows
\begin{align*}
{\rm I}_\tau \left(\vartheta \right)=\int_{\tau }^{T}\left[\frac{\dot
    f\left(\vartheta ,t\right) y_t\left(\vartheta \right)+f\left(\vartheta
    ,t\right)\dot y\left(\vartheta,t \right)}{\sigma \left(t\right)}\right]^2
{\rm d}t,\qquad {\rm I} \left(\vartheta \right)={\rm I}_0 \left(\vartheta \right).
\end{align*}

The family of measures of this statistical experiment is locally
asymptotically normal and therefore we have the lower bound on the mean
square risks of all estimators $\vartheta _\varepsilon $
\begin{align*}
\lim_{\nu \rightarrow 0}\Liminf_{\varepsilon \rightarrow
  0}\sup_{\left|\vartheta -\vartheta _0\right|<\nu } \varepsilon
^{-2}\Ex_\vartheta \left( \vartheta _\varepsilon-\vartheta  \right)^2\geq {\rm
I}\left(\vartheta _0\right)^{-1}.
\end{align*}
We call the estimator $\vartheta _\varepsilon ^*$ asymptotically efficient if
for all $\vartheta _0\in\Theta $ we have
\begin{align}
\label{14}
\lim_{\nu \rightarrow 0}\lim_{\varepsilon \rightarrow
  0}\sup_{\left|\vartheta -\vartheta _0\right|<\nu } \varepsilon
^{-2}\Ex_\vartheta \left( \vartheta _\varepsilon^*-\vartheta  \right)^2= {\rm
I}\left(\vartheta _0\right)^{-1}.
\end{align}

Introduce the estimator
\begin{align}
\label{15}
\vartheta _\varepsilon ^\star=\bar\vartheta_{\tau _\varepsilon }
+\frac{1}{{\rm I}_{\tau_\varepsilon }\left(\bar\vartheta_{\tau _\varepsilon } \right)}
\int_{\tau _\varepsilon }^{T}\frac{\dot M\left(\bar\vartheta_{\tau
    _\varepsilon },t \right) }{\sigma \left(t\right)^2}\left[{\rm
    d}X_t-f\left(\bar\vartheta_{\tau _\varepsilon
  },t\right)m\left(\bar\vartheta_{\tau _\varepsilon },t\right){\rm d}t\right] ,
\end{align}
where we denoted $\dot M\left(\vartheta,t \right)=\dot f\left(\vartheta
,t\right) y_t\left(\vartheta \right)+f\left(\vartheta ,t\right)\dot
y\left(\vartheta,t \right) $.

\begin{theorem}
\label{T2}
Suppose that the conditions ${\scr R}_1,\tilde{\scr R}_2,{\scr R}_3,{\scr
  R}_4$ be fulfilled and $\delta \in \left(0,1\right)$, then the One-step MLE
$\vartheta_{\varepsilon}^\star$ has the properties:
\begin{enumerate}
\item  It is consistent uniformly on compacts: for any $\nu >0$
$$ \sup_{\vartheta _0\in \KK}\Pb_{\vartheta_0}\left\{|\vartheta_{\varepsilon}^\star
-\vartheta_0|>\nu\right\}\longrightarrow0.
$$
\item It is uniformly  asymptotically normal
$$
  \varepsilon^{-1}\left(\vartheta_{\varepsilon}^\star-\vartheta_0\right)\Longrightarrow
  \zeta \sim \mathcal
N\left(0,{\rm I} \left(\vartheta_0 \right)^{-1}\right).
$$
\item The moments converge: for any $p>0$
\begin{align*}
\varepsilon ^{-p}\Ex_{\vartheta
  _0}\left|\vartheta_{\varepsilon}^\star-\vartheta_0  \right|^p\longrightarrow  \Ex_{\vartheta
  _0}\left|\zeta   \right|^p.
\end{align*}
\item It is  asymptotically efficient.
\end{enumerate}
\end{theorem}
{\bf Proof.} We have
\begin{align*}
&\frac{\vartheta_{\varepsilon}^\star-\vartheta_0}{\varepsilon }
= \frac{\bar\vartheta_{{\tau_\varepsilon}}-\vartheta_0 }{\varepsilon }
+{\rm I}_{\tau_\varepsilon }(\bar\vartheta_{\tau_\varepsilon})^{-1}\int_{\tau_\varepsilon}^T\frac{ \dot
  M\left(\bar\vartheta_{\tau _\varepsilon },t \right)
}{\sigma\left(t\right)}{\rm d}\bar W_t\\
&+\frac{1}{\varepsilon {\rm
  I}_{\tau_\varepsilon }(\bar\vartheta_{\tau_\varepsilon})}\int_{\tau_\varepsilon}^T\frac{ \dot
  M\left(\bar\vartheta_{\tau _\varepsilon },t \right)
}{\sigma\left(t\right)^2} \left[f(\vartheta_0,t)m(\vartheta_0,t)-
  f(\bar\vartheta_{\tau_\varepsilon},t)m(\bar\vartheta_{\tau_\varepsilon},t)\right]{\rm d}t.
\end{align*}
Here we used the {\it innovation representation} \cite{LS}
\begin{align*}
{\rm d}X_t=f\left(\vartheta _0\right)m\left(\vartheta _0,t\right){\rm
  d}t+\varepsilon \sigma \left(t\right) {\rm d}\bar W_t,\quad X_0=0,
\end{align*}
where the Wiener process $\bar W_t$ is defined by this equality.
Note that $\bar\vartheta_{\tau_\varepsilon}\rightarrow \vartheta_0$ and $\tau
_\varepsilon \rightarrow 0$, therefore uniformly on compacts
$$
\int_{\tau_\varepsilon}^T \frac{\dot M\left(\bar\vartheta_{\tau
    _\varepsilon },t \right)^2}{\sigma\left(t\right)^2}\;{\rm d}t\longrightarrow
\int_{0}^T \frac{\dot M\left(\vartheta_{0 },t
  \right)^2}{\sigma\left(t\right)^2}\;{\rm d}t= {\rm I}(\vartheta_{0})
$$
and by the central limit theorem we have
\begin{align*}
\Xi _\varepsilon =\frac{1}{ {\rm
  I}_{\tau_\varepsilon }(\bar\vartheta_{\tau_\varepsilon})}\int_{\tau_\varepsilon}^T
\frac{\dot M\left(\bar\vartheta_{\tau
    _\varepsilon },t \right)}{\sigma\left(t\right)}\;{\rm
  d}\bar W_t\Longrightarrow\mathcal N\left(0,{\rm
  I}(\vartheta_0)^{-1}\right).
\end{align*}
According to Taylor's formula
\begin{align*}
&f(\vartheta_0,t)m(\vartheta_0,t)-
f(\bar\vartheta_{\tau_\varepsilon},t)m(\bar\vartheta_{\tau_\varepsilon},t)\\
&\quad =\left[f(\vartheta_0,t)- f(\bar\vartheta_{\tau_\varepsilon},t)\right]
m(\vartheta_0,t)- f(\bar\vartheta_{\tau_\varepsilon},t)
\left[m(\bar\vartheta_{\tau_\varepsilon},t)- m(\vartheta_0,t)\right]\\
&\quad =\left(\vartheta_0-\bar\vartheta_{\tau_\varepsilon}\right)
\dot f(\tilde\vartheta_{\tau    _\varepsilon
},t)m(\vartheta_0,t)-\left(\bar\vartheta_{\tau_\varepsilon}-\vartheta_0\right)\dot
m(\tilde\vartheta_{\tau_\varepsilon},t  )f(\bar\vartheta_{\tau_\varepsilon},t)\\
&\quad =\left(\vartheta_0-\bar\vartheta_{\tau_\varepsilon}\right)
\left[\dot
  f(\bar\vartheta_{\tau_\varepsilon},t)m(\bar\vartheta_{\tau_\varepsilon},t
  )+f(\bar\vartheta_{\tau_\varepsilon},t)\dot
m(\bar\vartheta_{\tau_\varepsilon},t  ) + \bar R_\varepsilon\right]\\
&\quad =\left(\vartheta_0-\bar\vartheta_{\tau_\varepsilon}\right)
\left[\dot
  f(\bar\vartheta_{\tau_\varepsilon},t)y(\bar\vartheta_{\tau_\varepsilon},t
  )+f(\bar\vartheta_{\tau_\varepsilon},t)\dot
y(\bar\vartheta_{\tau_\varepsilon},t  ) + R_\varepsilon \right].
\end{align*}
Recall  that $\dot
m(\bar\vartheta_{\tau_\varepsilon},t  )=\dot
y(\bar\vartheta_{\tau_\varepsilon},t  )+O_p\left(\varepsilon \right) $. Here we denoted
$ \bar R_\varepsilon, R_\varepsilon,$
the random variables satisfying the estimates
\begin{align*}
\sup_{\vartheta _0\in \KK}\Ex_{\vartheta _0} \left| R_\varepsilon  \right|^p
\leq C\,\sup_{\vartheta _0\in \KK}\Ex_{\vartheta _0}
\left| \bar\vartheta_{\tau_\varepsilon}-\vartheta_0 \right|^p\leq C \,
\left(\frac{\varepsilon }{\sqrt{\tau _\varepsilon} }\right)^p
\end{align*}
Therefore we obtained the representation
\begin{align*}
\frac{\vartheta_{\varepsilon}^\star-\vartheta_0}{\varepsilon }
&=\Xi _\varepsilon +   \frac{\bar\vartheta_{{\tau_\varepsilon}}-\vartheta_0
}{\varepsilon } \left[1-\frac{1}{ {\rm
  I}_{\tau_\varepsilon }(\bar\vartheta_{\tau_\varepsilon})}\int_{\tau_\varepsilon}^T
\frac{\dot M\left(\bar\vartheta_{\tau
    _\varepsilon },t \right)^2}{\sigma\left(t\right)^2}\;{\rm
  d}t+R_\varepsilon \right]\\
&=\Xi _\varepsilon
+\frac{\left(\bar\vartheta_{{\tau_\varepsilon}}-\vartheta_0\right)R_\varepsilon
}{\varepsilon } ,
\end{align*}
where
\begin{align*}
\sup_{\vartheta _0\in \KK}\Ex_{\vartheta _0}\left| \frac{
  \left(\bar\vartheta_{{\tau_\varepsilon}}-\vartheta_0\right)R_\varepsilon
   }{\varepsilon }\right|^p&\leq C\,\sup_{\vartheta _0\in \KK} \Ex_{\vartheta _0}
\left|\frac{\left(\bar\vartheta_{{\tau_\varepsilon}}-\vartheta_0\ \right)^2}{\varepsilon
}\right|^p \leq
C\varepsilon ^{p\left(1-\delta\right) }\rightarrow 0.
\end{align*}
Hence
\begin{align*}
\frac{\vartheta_{\varepsilon}^\star-\vartheta_0}{\varepsilon }
&=\Xi _\varepsilon +o\left(1\right)\Longrightarrow {\cal N}\left(0, {\rm I}\left(\vartheta _0\right)^{-1}\right).
\end{align*}
From this proof follows the uniform on compacts  convergence of moments too
\begin{align*}
{\varepsilon }^{-2}\Ex_{\vartheta _0}
\left(\vartheta_{\varepsilon}^\star-\vartheta_0\right)^2\longrightarrow {\rm
    I}\left(\vartheta _0\right)^{-1}.
\end{align*}
This uniform convergence allows us to write
\begin{align*}
\sup_{\left|\vartheta -\vartheta _0\right|<\nu }{\varepsilon }^{-2}\Ex_{\vartheta}
\left(\vartheta_{\varepsilon}^\star-\vartheta\right)^2\longrightarrow \sup_{\left|\vartheta -\vartheta _0\right|<\nu }{\rm
    I}\left(\vartheta\right)^{-1}\longrightarrow {\rm
    I}\left(\vartheta_0\right)^{-1},
\end{align*}
where the last limit was obtained as $\nu\rightarrow 0 $.  Hence we verified
the asymptotic efficiency \eqref{14} of the One-step MLE $\vartheta _\varepsilon ^\star$.

\bigskip

 {\bf Case} $f\left(\vartheta ,t\right)=f\left(t\right)$.

The estimator $\bar\vartheta _{\tau _\varepsilon }$ is defined by the same
equation  \eqref{8}.
The Fisher information is
\begin{align*}
{\rm I}_{\tau  }\left(\vartheta \right)=\int_{\tau }^{T}\frac{f\left(t\right)^2\dot y\left(\vartheta
  ,t\right)^2}{ \sigma \left(t\right)^2}{\rm d}t,\qquad {\rm I}
\left(\vartheta \right)={\rm I}_{0 }\left(\vartheta \right)
\end{align*}
Introduce the estimator
\begin{align}
\label{16}
\vartheta _\varepsilon ^\star=\bar\vartheta_{\tau _\varepsilon }
+\frac{1}{{\rm I}_{\tau_\varepsilon }\left(\bar\vartheta_{\tau _\varepsilon } \right)}
\int_{\tau _\varepsilon }^{T}\frac{f\left(t\right)\dot y\left(\bar\vartheta_{\tau
    _\varepsilon },t \right) }{\sigma \left(t\right)^2}\left[{\rm
    d}X_t-f\left(t\right)m\left(\bar\vartheta_{\tau _\varepsilon },t\right){\rm d}t\right]
\end{align}
and condition

${\scr R}_2^*.$ {\it The function $
  a\left(\vartheta ,\cdot  \right) \in {\cal C}^{\left(2\right)}_\vartheta $. The function $\dot
  a\left(\vartheta ,t  \right)\in {\cal C}^{\left(1\right)}_t$.}

\begin{proposition}
\label{P2}
Suppose that the conditions ${\scr R}_1,{\scr R}_2^*,{\scr R}_3',{\scr
  R}_4$ be fulfilled and $\delta \in \left(0,1/3\right)$, then the One-step MLE
$\vartheta_{\varepsilon}^\star$  has the properties
\begin{enumerate}
\item It is  uniformly on compacts consistent: for any $\nu >0$
$$ \sup_{\vartheta _0\in \KK}\Pb_{\vartheta_0}\left\{|\vartheta_{\varepsilon}^\star
-\vartheta_0|>\nu\right\}\longrightarrow0.
$$
\item It is uniformly  asymptotically normal
$$
\varepsilon^{-1}\left(\vartheta_{\varepsilon}^\star-\vartheta_0\right)\Longrightarrow
\mathcal
N\left(0,{\rm I} \left(\vartheta_0 \right)^{-1}\right),
$$
\item It is asymptotically efficient.
\end{enumerate}
\end{proposition}
{\bf Proof.} The proof of this proposition is close to the given above proof
of Theorem \ref{T1}. The difference is in the estimation of the term
\begin{align*}
\Ex_{\vartheta _0} \frac{\left|
  \left(\bar\vartheta_{{\tau_\varepsilon}}-\vartheta_0\right)R_\varepsilon
  \right| }{\varepsilon }\leq C\varepsilon ^{-1} \Ex_{\vartheta _0}
\left(\bar\vartheta_{{\tau_\varepsilon}}-\vartheta_0\ \right)^2 \leq
C\varepsilon ^{-1+2-3\delta } =C\varepsilon ^{1-3\delta }\rightarrow 0.
\end{align*}

\section{One-step MLE-process}

Let us consider slightly different problem. We have the model of observations
\eqref{1}-\eqref{2} with unknown parameter $\vartheta \in\Theta $ and we are
interested in the construction of the adaptive Kalman-Bucy filter to
approximate $m\left(\vartheta _0,t\right)=\Ex_{\vartheta
  _0}\left(Y_t|X_s,0\leq s\leq t\right)$. Of course, we can not use
$m\left(\vartheta _\varepsilon ^\star,t\right) $ because the estimator
$\vartheta _\varepsilon ^\star$ depends on the observations
$X^T=\left(X_s,0\leq s\leq T\right)$, where the part $\left(X_s,t< s\leq T\right)$
are from the future. Therefore we need an estimator-process $\vartheta
_{t,\varepsilon} ^\star, 0<t\leq T$, where the estimator $\vartheta
_{t,\varepsilon} ^\star $ has the following properties
\begin{itemize}
\item {\it  it depends on $\left(X_s,\, 0\leq s\leq t\right)$;
\item  it can be  easily  calculated;
 \item  it is asymptotically efficient.}
\end{itemize}
We can not use the MLE $\hat \vartheta _{t,\varepsilon} $ obtained as  solution of the
equation
\begin{align*}
L(\hat \vartheta _{t,\varepsilon} ,X^t )=\sup_{\vartheta \in\Theta
}L\left(  \vartheta   ,X^t \right)
\end{align*}
  and put
$m(\hat\vartheta _{t,\varepsilon },t) $ because its
calculation for all $t\in\left(0,T\right)$ requires calculation of the solutions
of the filtration equations \eqref{5}-\eqref{6}. Of course numerical
realization of such procedure is too complicated.

We propose a modified One-step MLE $\vartheta _{t,\varepsilon }^\star,\tau
_\varepsilon \leq t\leq T$ called {\it One-step MLE-process} defined as follows
\begin{align}
\label{17}
\vartheta _{t,\varepsilon }^\star=\bar\vartheta _{\tau _\varepsilon
}+\frac{1}{{\rm I}_{\tau _\varepsilon }^t\left(\bar\vartheta _{\tau _\varepsilon
}\right)}\int_{\tau _\varepsilon }^{t}\frac{\dot M\left(\bar\vartheta _{\tau _\varepsilon
},s \right)}{\sigma \left(s\right)^2}\left[{\rm d}X_s-f\left(\bar\vartheta
  _{\tau _\varepsilon }, s\right) m\left(\bar\vartheta
  _{\tau _\varepsilon }, s\right){\rm d}s\right].
\end{align}
Here
\begin{align*}
{\rm I}_{\tau }^t\left(\vartheta\right)=\int_{\tau }^{t}\frac{\dot
  M\left(\vartheta ,s\right)^2}{\sigma \left(s\right)^2} \;{\rm d}s,\qquad
{\rm I}^t\left(\vartheta\right)={\rm I}_{0 }^t\left(\vartheta\right) .
\end{align*}
\begin{proposition}
\label{P3} Let the conditions  ${\scr R}_1,{\scr R}_2^*,{\scr R}_3',{\scr
  R}_4$ be fulfilled. Then the One-step MLE-process
$\vartheta _{t,\varepsilon }^\star,\tau _\varepsilon <t\leq T $ is consistent,
asymptotically normal
\begin{align*}
\frac{\vartheta _{t,\varepsilon }^\star-\vartheta _0}{\varepsilon
}\Longrightarrow {\cal N}\left(0,{\rm I}_{0}^t\left(\vartheta_0\right)^{-1} \right)
\end{align*}
and the random process $\eta _\varepsilon \left(t\right)=\varepsilon
^{-1}\left( \vartheta _{t,\varepsilon }^\star-\vartheta _0\right), \tau \leq
t\leq T$ for any $\tau \in (0,T)$ converges in distribution to the random
process
\begin{align*}
\eta _\varepsilon \left(t\right)\Longrightarrow \eta \left(t\right)={\rm
  I}_{0}^t\left(\vartheta_0\right)^{-1}\int_{0}^{t}\frac{\dot M\left(\vartheta
  _0,s\right)}{\sigma \left(s\right)}{\rm d}w_s,\qquad \tau \leq t\leq T,
\end{align*}
where $w_s,0\leq s\leq T$ is some Wiener process.
\end{proposition}
{\bf Proof.} As $\tau _\varepsilon \rightarrow 0$ we have $\tau _\varepsilon
<\tau $ for any $\tau >0$ starting from the corresponding value, say,
$\varepsilon _0$. For any $t\in (\tau _\varepsilon ,T]$ according to Theorem
  \ref{T1} (under corresponding conditions) this estimator is consistent and
  asymptotically normal. Moreover, if we consider the vector $\eta
  _\varepsilon \left(t_1\right),\ldots , \eta _\varepsilon \left(t_k\right) $,
  then repeating the proof of Theorem \ref{T1} we obtain the convergence of
  its distribution to the distribution of the vector $\eta
  \left(t_1\right),\ldots , \eta \left(t_k\right) $. Further, the obtained in
  the proof of Theorem \ref{T1} representation
\begin{align*}
&\eta _\varepsilon \left(t\right)
= \frac{\bar\vartheta_{{\tau_\varepsilon}}-\vartheta_0 }{\varepsilon }
+{\rm I}_{\tau_\varepsilon }^t(\bar\vartheta_{\tau_\varepsilon})^{-1}\int_{\tau_\varepsilon}^t\frac{ \dot
  M\left(\bar\vartheta_{\tau _\varepsilon },s \right)
}{\sigma\left(s\right)}{\rm d}\bar W_s\\
&+\frac{1}{\varepsilon {\rm
  I}_{\tau_\varepsilon }^t(\bar\vartheta_{\tau_\varepsilon})}\int_{\tau_\varepsilon}^t\frac{ \dot
  M\left(\bar\vartheta_{\tau _\varepsilon },s \right)
}{\sigma\left(s\right)^2} \left[f(\vartheta_0,s)m(\vartheta_0,s)-
  f(\bar\vartheta_{\tau_\varepsilon},s)m(\bar\vartheta_{\tau_\varepsilon},s)\right]{\rm d}s.
\end{align*}
 allows us to verify the estimate
\begin{align}
\label{18}
\Ex_{\vartheta _0}\left|\eta _\varepsilon \left(t_1\right)-\eta _\varepsilon
\left(t_2\right)   \right|^2\leq C\left|t_1-t_2\right|^2 ,
\end{align}
where the constant $C>0$ does not depend on $\varepsilon $.  The calculations
are direct but cumbersome. The convergence of
finite-dimensional distributions and estimate \eqref{18} provide us the weak
convergence of the measures induced in the space of continuous functions
${\cal C}\left[\tau ,T\right]$ by the processes $\eta
_\varepsilon \left(\cdot \right)$ to the measure of the process $ \eta
\left(\cdot \right)$ (see details of the proof in the similar situation in
\cite{Kut17}).

This weak convergence provides us the relation: for any $\tau \in
\left(0,T\right)$ and any $\nu >0$
\begin{align*}
&\Pb_{\vartheta _0}\left(\sup_{\tau \leq t\leq T}{\rm I}_{\tau_\varepsilon
  }^t(\bar\vartheta_{\tau_\varepsilon}) \left|\frac{\vartheta _{t,\varepsilon
    }^\star-\vartheta _0}{\varepsilon } \right|>\nu \right)\\ &\qquad \qquad
  \longrightarrow \Pb_{\vartheta _0}\left(\sup_{\tau \leq t\leq
    T}\left|\int_{0}^{t}\frac{\dot M\left(\vartheta _0,s\right)}{\sigma
    \left(s\right)}{\rm d}w_s \right|>\nu \right)\\ &\qquad \qquad
  = \Pb_{\vartheta _0}\left(\sup_{\Lambda_\tau \leq \lambda \leq
   \Lambda_T}\left|W\left(\lambda \right) \right|>\nu \right).
\end{align*}
Here $W\left(\cdot \right)$ is a Wiener process and
\begin{align*}
\lambda=\int_{0}^{t}\frac{\dot M\left(\vartheta
  _0,s\right)^2}{\sigma \left(s\right)^2} \;{\rm d}s,\quad \Lambda_\tau
=\int_{0}^{\tau }\frac{\dot M\left(\vartheta
  _0,s\right)^2}{\sigma \left(s\right)^2} \;{\rm d}s,\quad \Lambda_T
=\int_{0}^{T }\frac{\dot M\left(\vartheta
  _0,s\right)^2}{\sigma \left(s\right)^2} \;{\rm d}s.
\end{align*}

Note that the One-step MLE-process is uniformly consistent: for any $\tau
\in\left(0,T\right)$ and any $\nu >0$ we have
\begin{align*}
\Pb_{\vartheta _0}\left(\sup_{\tau \leq t\leq T} \left|\vartheta _{t,\varepsilon
    }^\star-\vartheta _0 \right|>\nu \right)\longrightarrow 0.
\end{align*}

\bigskip

The adaptive equations of filtration can be written as follows
\begin{align}
{\rm d}m_\varepsilon^\star \left(t\right)&=a\left(\vartheta _{t,\varepsilon }^\star
  ,t\right)m_\varepsilon ^\star\left(t\right){\rm d}t\nonumber\\
&\qquad \qquad +\frac{\gamma_{*,\varepsilon }
  \left(t\right)f\left(\vartheta _{t,\varepsilon }^\star
  ,t\right) }{\sigma \left(t\right)^2} \left[{\rm
  d}X_t-f\left(\vartheta _{t,\varepsilon }^\star ,t\right)m_\varepsilon^\star
  \left(t\right){\rm d}t \right],\label{19}\\
 \frac{\partial \gamma_{*,\varepsilon } \left(t\right)}{\partial
  t}&=2a\left(\vartheta _{t,\varepsilon }^\star ,t\right)\gamma
 _{*,\varepsilon } \left(t\right)-\frac{\gamma _{*,\varepsilon }
  \left(t\right)^2f\left(\vartheta _{t,\varepsilon }^\star
  ,t\right)^2}{\sigma \left(t\right)^2}+b\left(t\right)^2,
\label{20}
\end{align}
with initial values $m_\varepsilon \left(\tau _\varepsilon \right)=
m_\varepsilon \left(\vartheta ,\tau _\varepsilon \right)$ and
$\gamma_{*,\varepsilon } \left(t\right)=0 $.

 The equations \eqref{18}-\eqref{20} give us closed system, which allows us to
 calculate the estimator and the approximation of the conditional expectation in
 recurrent form.

\section{On  Efficient Estimation of $m\left(\vartheta
_0,t\right)$}

Recall that $m\left(\vartheta _0,t\right)$ is mean squared optimal estimator
of $Y_t$. The random process $m_\varepsilon \left(t\right),\tau _\varepsilon
\leq t\leq T$ can be considered as estimator of the random function
$m\left(\vartheta _0,t\right),\tau _\varepsilon \leq t\leq T$ and it is
interesting to see what are asymptotically efficient estimators in this
problem.

We have the following lower bound on the mean square error of any estimator  $\bar m_\varepsilon
\left(t\right)$    of the random function
$m\left(\vartheta ,t\right)$.

\begin{theorem}
\label{T3}
 Let the conditions  ${\scr R}_1,\tilde{\scr R}_2,{\scr R}_3,{\scr
  R}_4$ be fulfilled and $\delta \in \left(0,1\right)$ be fulfilled, then for any $\vartheta _0
\in \Theta $ and any estimator  $\bar m_\varepsilon \left(t\right)$ we have
\begin{align}
\label{21}
\lim_{\nu \rightarrow 0}\Liminf_{\varepsilon \rightarrow
  0}\sup_{\left|\vartheta -\vartheta _0\right|<\nu }\varepsilon ^{-2} \Ex_\vartheta \left|\bar m_\varepsilon
\left(t\right)-m\left(\vartheta ,t\right) \right|^2\geq \frac{\dot
  y\left(\vartheta _0,t\right)^2}{ {\rm I}^t\left(\vartheta _0\right)}.
\end{align}
\end{theorem}
{\bf Proof.} The proof of this inequality follows the main steps of the proof
of van Trees inequality given in \cite{GL95} and almost coincides with the
proof of similar bound in \cite{KZ14}. Let us remind here some steps of the
proof. Introduce a  density function $p\left(\vartheta \right),\vartheta
_0-\nu <\vartheta <\vartheta _0+\nu$ such that  $p\left(\vartheta_0\pm \nu
\right)=0 $ and  Fisher information
\begin{align*}
{\rm I}_p=\int_{\vartheta _0-\nu }^{\vartheta _0+\nu}\frac{\dot
  p\left(\vartheta \right)^2}{p\left(\vartheta \right)}\;{\rm
  d}\vartheta<\infty   .
\end{align*}
  Then we can write
\begin{align*}
\sup_{\left|\vartheta -\vartheta _0\right|<\nu }
\Ex_\vartheta \left|\bar m_\varepsilon \left(t\right)-m\left(\vartheta
,t\right) \right|^2&\geq \int_{\vartheta _0-\nu }^{\vartheta _0+\nu
}\Ex_\vartheta \left|\bar m_\varepsilon \left(t\right)-m\left(\vartheta
,t\right) \right|^2p\left(\vartheta \right){\rm d}\vartheta\\
&=\EE \left|\bar m_\varepsilon \left(t\right)-m\left(\vartheta
,t\right) \right|^2,
\end{align*}
where we denote $\EE$ double mathematical expectation defined by the last
equality. Further we need a version of van Trees inequality \cite{GL95}, which is obtained
as follows. Let us denote
\begin{align*}
L\left(\vartheta ,\vartheta
_0,X^t\right)&=\exp\left\{\int_{0}^{t}\frac{f\left(\vartheta ,s\right)m\left(\vartheta
  ,s\right)-f\left(\vartheta_0 ,s\right)m\left(\vartheta_0 ,s\right)}{\varepsilon ^2\sigma
  \left(s\right)^2}\;{\rm d}X_s\right.\\
&\qquad \quad \left. - \int_{0}^{t}\frac{f\left(\vartheta ,s\right)^2m\left(\vartheta
  ,s\right)^2-f\left(\vartheta_0 ,s\right)^2m\left(\vartheta_0 ,s\right)^2}{2\varepsilon ^2\sigma
  \left(s\right)^2}\;{\rm d}s \right\} .
\end{align*}
We have
\begin{align*}
&\int_{\vartheta _0-\nu }^{\vartheta _0+\nu}\left[\bar m_\varepsilon \left(t\right)-m\left(\vartheta ,s\right) \right]\frac{\partial }{\partial \vartheta }\left[ L\left(\vartheta ,\vartheta
_0,X^t\right)   p\left(\vartheta \right)\right]{\rm d}\vartheta \\
&\qquad =\left(\bar m_\varepsilon \left(t\right)-m\left(\vartheta ,s\right) \right)\left.L\left(\vartheta ,\vartheta
_0,X^t\right)   p\left(\vartheta \right)\right|_{\vartheta _0-\nu }^{\vartheta
    _0+\nu }\\
&\qquad\qquad\qquad+\int_{\vartheta _0-\nu }^{\vartheta _0+\nu} \dot
  m\left(\vartheta ,s\right)   L\left(\vartheta ,\vartheta
_0,X^t\right)   p\left(\vartheta \right) {\rm d}\vartheta\\
&\qquad =\int_{\vartheta _0-\nu }^{\vartheta _0+\nu} \dot
  m\left(\vartheta ,s\right)   L\left(\vartheta ,\vartheta
_0,X^t\right)   p\left(\vartheta \right) {\rm d}\vartheta.
\end{align*}
Therefore
\begin{align}
&\Ex_{\vartheta _0}\int_{\vartheta _0-\nu }^{\vartheta _0+\nu}\left[\bar
    m_\varepsilon \left(t\right)-m\left(\vartheta ,s\right)
    \right]\frac{\partial }{\partial \vartheta }\left[ L\left(\vartheta
    ,\vartheta _0,X^t\right) p\left(\vartheta \right)\right]{\rm
    d}\vartheta\nonumber\\
 &\qquad \qquad =\Ex_{\vartheta _0}\int_{\vartheta _0-\nu
  }^{\vartheta _0+\nu} \dot m\left(\vartheta ,s\right) L\left(\vartheta
  ,\vartheta _0,X^t\right) p\left(\vartheta \right){\rm d}\vartheta\nonumber\\
 &\qquad \qquad =\int_{\vartheta _0-\nu
  }^{\vartheta _0+\nu} \Ex_{\vartheta}\,\dot m\left(\vartheta ,s\right)
  p\left(\vartheta \right){\rm d}\vartheta,
\label{22}
\end{align}
where we changed the measure $\Ex_{\vartheta _0}L\left(\vartheta
  ,\vartheta _0,X^t\right)=\Ex_{\vartheta} $. Recall that
\begin{align*}
&\frac{\partial }{\partial \vartheta }\ln L\left(\vartheta  ,\vartheta
_0,X^t\right)\\
&\qquad =\int_{0}^{t}\frac{f\left(\vartheta ,s\right) \dot m\left(\vartheta
  ,s\right)+\dot f\left(\vartheta ,s\right)  m\left(\vartheta
  ,s\right)}{\varepsilon^2\sigma \left(s\right) ^2
}\left[{\rm d}X_s- f\left(\vartheta ,s\right)m\left(\vartheta ,s\right){\rm d}s\right].
\end{align*}
Hence (below $L_{\vartheta ,\vartheta _0}p_\vartheta =L\left(\vartheta
  ,\vartheta _0,X^t\right) p\left(\vartheta \right)$)
\begin{align}
&\left(\Ex_{\vartheta _0} \int_{\vartheta _0-\nu }^{\vartheta
    _0+\nu}\left[\bar m_\varepsilon \left(t\right)-m\left(\vartheta ,t\right)
    \right]\frac{\partial }{\partial \vartheta }\left[ L\left(\vartheta
    ,\vartheta _0,X^t\right) p\left(\vartheta \right)\right]{\rm d}\vartheta
  \right)^2\nonumber\\ 
&\qquad =\left(\Ex_{\vartheta _0} \int_{\vartheta _0-\nu
  }^{\vartheta _0+\nu}\left[\bar m_\varepsilon \left(t\right)-m\left(\vartheta
    ,t\right) \right]\frac{\partial \ln \left[L_{\vartheta ,\vartheta
        _0}p_\vartheta\right] }{\partial \vartheta } L_{\vartheta ,\vartheta
    _0} p_\vartheta \,{\rm d}\vartheta \right)^2\nonumber\\ 
&\qquad\leq \Ex_{\vartheta
    _0} \int_{\vartheta _0-\nu }^{\vartheta _0+\nu}\left[\bar m_\varepsilon
    \left(t\right)-m\left(\vartheta ,t\right) \right]^2 L_{\vartheta
    ,\vartheta _0}p_\vartheta \, {\rm d}\vartheta \nonumber\\ 
&\qquad \qquad \qquad
  \quad \times\Ex_{\vartheta _0} \int_{\vartheta _0-\nu }^{\vartheta _0+\nu}
  \left[\frac{\partial \ln \left[L_{\vartheta ,\vartheta _0}p_\vartheta\right]
    }{\partial \vartheta } \right] ^2 L_{\vartheta ,\vartheta _0}p_\vartheta
       {\rm d}\vartheta\nonumber\\ 
&\qquad = \int_{\vartheta _0-\nu }^{\vartheta
         _0+\nu}\Ex_{\vartheta}\left[\bar m_\varepsilon
         \left(t\right)-m\left(\vartheta ,t\right) \right]^2
       p\left(\vartheta\right) \, {\rm d}\vartheta\nonumber\\
&\qquad \qquad \qquad \qquad \times\left( \int_{\vartheta _0-\nu
       }^{\vartheta _0+\nu}\Ex_{\vartheta } \left[\frac{\partial \ln
           L\left(\vartheta ,\vartheta _0 , X^t\right) }{\partial \vartheta } \right] ^2
        p\left(\vartheta\right) \, {\rm d}\vartheta+{\rm I}_p \right).
\label{23}
\end{align}
Here we used Cauchy-Shwarz inequality and the property of stochastic integral
\begin{align*}
&\Ex_{\vartheta } \frac{\partial \ln
           L\left({\vartheta ,\vartheta _0},X^t\right) }{\partial \vartheta}\\
&\;=\Ex_{\vartheta } \int_{0}^{t }\frac{\dot f\left(\vartheta
  ,s\right)m\left(\vartheta ,s\right)+f\left(\vartheta
  ,s\right)\dot m\left(\vartheta ,s\right) }{\varepsilon ^2\sigma
  \left(s\right)^2}\left[{\rm d}X_s-f\left(\vartheta
  ,s\right)m\left(\vartheta ,s\right){\rm d}s \right]= 0.
\end{align*}
We have
\begin{align}
&\Ex_{\vartheta } \left[\frac{\partial \ln
         L\left({\vartheta ,\vartheta _0},X^t\right) }{\partial \vartheta }
  \right] ^2\nonumber\\
&\quad =\Ex_{\vartheta } \left(\int_{0}^{t }\frac{\dot f\left(\vartheta
  ,s\right)m\left(\vartheta ,s\right)+f\left(\vartheta
  ,s\right)\dot m\left(\vartheta ,s\right) }{\varepsilon\sigma
  \left(s\right)}{\rm d}\bar W_s\right)^2\nonumber\\
&\quad =\Ex_{\vartheta } \int_{0}^{t }\left(\frac{\dot f\left(\vartheta
  ,s\right)m\left(\vartheta ,s\right)+f\left(\vartheta
  ,s\right)\dot m\left(\vartheta ,s\right) }{\varepsilon\sigma
  \left(s\right)}\right)^2{\rm d}s\nonumber\\
&\quad = \frac{1}{\varepsilon ^2}\int_{0}^{t }\left(\frac{\dot f\left(\vartheta
  ,s\right)y_s\left(\vartheta \right)+f\left(\vartheta
  ,s\right)\dot y\left(\vartheta ,s\right) }{\sigma
  \left(s\right)}\right)^2{\rm d}s\left(1+O\left(\varepsilon \right)\right)\nonumber\\
&\quad = \frac{1}{\varepsilon ^2}\int_{0}^{t }\left(\frac{\dot M\left(\vartheta
  ,s\right) }{\sigma
  \left(s\right)}\right)^2{\rm d}s\left(1+O\left(\varepsilon \right)\right)=
  \frac{1}{\varepsilon ^2} {\rm I}^t\left(\vartheta
  \right)\left(1+O\left(\varepsilon \right)\right).
\label{24}
\end{align}
Using \eqref{22}-\eqref{24} we obtain the relation
\begin{align*}
&\left(\int_{\vartheta _0-\nu
  }^{\vartheta _0+\nu} \Ex_{\vartheta}\,\dot m\left(\vartheta ,s\right)
  p\left(\vartheta \right){\rm d}\vartheta \right)^2\\
&\qquad \qquad \leq \EE \left|\bar
  m_\varepsilon \left(t\right)-m\left(\vartheta ,t\right)\right|^2
  \left(\int_{\vartheta_0-\nu}^{\vartheta_0+\nu} \frac{{\rm I}^t\left(\vartheta
  \right)}{\varepsilon ^2}\left(1+O\left(\varepsilon \right)\right)p(\vartheta){\rm d}\vartheta+{\rm I}_p \right),
\end{align*}
which can be written as follows
\begin{align*}
&\varepsilon ^{-2}\EE \left|\bar m_\varepsilon \left(t\right)-m\left(\vartheta
  ,t\right)\right|^2\geq \frac{ \left(\int_{\vartheta _0-\nu }^{\vartheta
      _0+\nu} \dot y\left(\vartheta ,s\right) p\left(\vartheta \right){\rm
      d}\vartheta \right)^{2}\left(1+O\left(\varepsilon
    \right)\right)}{\int_{\vartheta _0-\nu }^{\vartheta _0+\nu}{\rm
      I}^t\left(\vartheta \right)p\left(\vartheta \right){\rm d}\vartheta
    \left(1+O\left(\varepsilon \right)\right)+\varepsilon ^{2}{\rm I}_p}
  \\
&\qquad \qquad \longrightarrow \frac{ \left(\int_{\vartheta _0-\nu
    }^{\vartheta _0+\nu} \dot y\left(\vartheta ,s\right) p\left(\vartheta
    \right){\rm d}\vartheta \right)^{2}}{\int_{\vartheta _0-\nu }^{\vartheta
      _0+\nu}{\rm I}^t\left(\vartheta \right)p\left(\vartheta \right){\rm
      d}\vartheta }\longrightarrow \frac{ \dot y\left(\vartheta_0 ,t\right)
    ^{2}}{{\rm I}^t\left(\vartheta _0 \right) }
\end{align*}
where the first and second limits correspond  to $\varepsilon \rightarrow 0$
and $\nu \rightarrow 0$.  Therefore we obtained \eqref{21}.

This bound allows us to define asymptotically efficient estimator of the
conditional expectation
$m\left(\vartheta _0,t\right)$ as  estimator $\hat m_\varepsilon
\left(t\right)$ satisfying
\begin{align}
\label{25}
\lim_{\nu \rightarrow 0}\lim_{\varepsilon \rightarrow
  0}\sup_{\left|\vartheta -\vartheta _0\right|<\nu }\varepsilon ^{-2}
\Ex_\vartheta \left|\hat m_\varepsilon
\left(t\right)-m\left(\vartheta ,t\right) \right|^2= \frac{\dot
  y\left(\vartheta _0,t\right)^2}{ {\rm I}^t\left(\vartheta _0\right)}
\end{align}
 for all $\vartheta _0\in \Theta $.

To construct asymptotically efficient estimator we slightly modify the
estimator $m_\varepsilon \left(t\right) $.

The solutions of the equations \eqref{5} and  \eqref{19} can be written as
follows
\begin{align}
\label{26}
 m\left(\vartheta _0,t\right)&=y_0\,N
 \left(\vartheta_0,t\right)+N
 \left(\vartheta_0,t\right)\int_{\tau _\varepsilon
 }^{t}\frac{Q \left(\vartheta_0,s\right) }{N
   \left(\vartheta_0 ,s\right)}\;{\rm d}X_s,\\
 m_\varepsilon \left(t\right)&=y_0\,N_\varepsilon
 \left(\vartheta^\star,t\right)+N_\varepsilon
 \left(\vartheta^\star,t\right)\int_{\tau _\varepsilon
 }^{t}\frac{Q_\varepsilon \left(\vartheta^\star,s\right) }{N_\varepsilon
   \left(\vartheta^\star ,s\right)}\;{\rm d}X_s,\nonumber
\end{align}
where
\begin{align*}
N \left(\vartheta_0,t\right)&=\exp\left\{\int_{\tau _\varepsilon
}^{t}\left[a\left(\vartheta _0,s\right)-{\gamma
  \left(\vartheta _0,s\right)f\left(\vartheta _0
  ,s\right)^2 }{\sigma \left(s\right)^{-2}} \right] {\rm d}s \right\},\\
N_\varepsilon \left(\vartheta^\star,t\right)&=\exp\left\{\int_{\tau _\varepsilon
}^{t}\left[a\left(\vartheta _{s,\varepsilon }^\star,s\right)-{\gamma_{*,\varepsilon }
  \left(s\right)f\left(\vartheta _{s,\varepsilon }^\star
  ,s\right)^2 }{\sigma \left(s\right)^{-2}} \right] {\rm d}s \right\},\\
Q \left(\vartheta_0,t\right)&={\gamma
  \left(\vartheta _0,t\right)f\left(\vartheta _0
  ,t\right) }{\sigma \left(t\right)^{-2}},\quad Q_\varepsilon
\left(\vartheta^\star,t\right)={\gamma_{*,\varepsilon }
  \left(t\right)f\left(\vartheta _{t,\varepsilon }^\star
  ,t\right) }{\sigma \left(t\right)^{-2}}.
\end{align*}
We can consider another possibility of the approximation of $m\left(\vartheta
_0,t\right)$ using the One-step MLE and equation \eqref{19}. We can not put
$\vartheta _{t,\varepsilon }^\star$ in  $m\left(\vartheta
_0,t\right) $ because the stochastic integral
\begin{align*}
\int_{\tau _\varepsilon
 }^{t}F\left(\vartheta _{t,\varepsilon }^\star,s\right) \;{\rm d}X_s
\end{align*}
is not well defined for $F\left(\vartheta ,t\right)={Q \left(\vartheta,s\right) }{N
   \left(\vartheta ,s\right)^{-1}}$.  We will use the so-called robust version of the
integral. It is given by the formula
\begin{align}
\label{27}
\int_{\tau _\varepsilon
 }^{t}F\left(\vartheta ,s\right)\;{\rm d}X_s=F\left(\vartheta
,t\right)X_t-F\left(\vartheta ,\tau _\varepsilon\right)X_{\tau _\varepsilon }-\int_{\tau _\varepsilon
 }^{t}F'\left(\vartheta ,s\right)X_s\;{\rm d}s,
\end{align}
where
\begin{align}
\label{28}
F'\left(\vartheta ,s\right)&=\frac{Q'\left(\vartheta ,s\right)N\left(\vartheta
  ,s\right)-Q\left(\vartheta ,s\right)N'\left(\vartheta ,s\right) }{N\left(\vartheta
  ,s\right)^2},\\
 N'\left(\vartheta ,s\right)&=N\left(\vartheta
,s\right)\left[a\left(\vartheta ,s\right)-\gamma \left(\vartheta
  _0,s\right)f\left(\vartheta _0,s\right)^2\sigma \left(s\right)^{-2}\right] \nonumber.
\end{align}

Let us denote the the of this equality as $G\left(\vartheta ,X,t\right)$ and
introduce the estimator
\begin{align}
\label{29}
m_\varepsilon^\star \left(t\right)=y_0\,N\left(\vartheta _{t,\varepsilon
}^\star,t\right)+N\left(\vartheta _{t,\varepsilon }^\star,t\right)
G\left(\vartheta _{t,\varepsilon }^\star ,X,t\right) ,
\end{align}
which we compare with
\begin{align}
\label{30}
m \left(\vartheta _0,t\right)=y_0\,N\left(\vartheta
_0,t\right)+N\left(\vartheta_0,t\right)
G\left(\vartheta_0 ,X,t\right) .
\end{align}
Remind that we comprehend $F'\left(\vartheta_{t,\varepsilon }^\star ,s\right)$ as that is written in \eqref{28}, therefore the equality
\eqref{27} will no more be valid if $\vartheta $ is replaced by $\vartheta_{t,\varepsilon
}^\star $. This means that the estimator in \eqref{29} is not
exactly $m\left(\vartheta_{t,\varepsilon }^\star,t \right)$.

 Recall that
\begin{align*}
\left.m \left(\vartheta _0,t\right)\right|_{\varepsilon =0}&=
y_0\,N\left(\vartheta _0,t\right)+N\left(\vartheta_0,t\right)
G\left(\vartheta_0 ,x,t\right)=y_t\left(\vartheta _0\right).
\end{align*}

It is interesting to see if the proposed estimator $m_\varepsilon^\star
\left(t\right)$ is optimal in some sense.

\begin{proposition}
\label{P4} Let the conditions  ${\scr R}_1,{\scr R}_2^*,{\scr R}_3',{\scr
  R}_4$ be fulfilled then the estimator $m_\varepsilon^\star
\left(t\right)$ is asymptotically efficient, i.e., for any $t\in (0,T]$ and
  any $\vartheta
_0\in\Theta $ we have
 \begin{align*}
\lim_{\nu \rightarrow 0}\lim_{\varepsilon \rightarrow 0}\sup_{\left|\vartheta
  -\vartheta _0\right|<\nu }\varepsilon ^{-2} \Ex_\vartheta \left|
m_\varepsilon^\star \left(t\right)-m\left(\vartheta ,t\right) \right|^2=
\frac{\dot y\left(\vartheta _0,t\right)^2}{ {\rm I}^t\left(\vartheta
  _0\right)}.
\end{align*}
\end{proposition}
{\bf Proof.} We can write formally
\begin{align*}
m_\varepsilon^\star \left(t\right)-m\left(\vartheta ,t\right)&=m
\left(\vartheta _{t,\varepsilon }^\star, t\right)-m\left(\vartheta
,t\right)=\left(\vartheta _{t,\varepsilon }^\star-\vartheta  \right)\dot m\left(\vartheta
,t\right) \left(1+o\left(1\right)\right)\\
&= \left(\vartheta _{t,\varepsilon }^\star-\vartheta  \right)\dot y\left(\vartheta
,t\right) \left(1+o\left(1\right)\right)
\end{align*}
and
\begin{align*}
&\sup_{\left|\vartheta -\vartheta _0\right|<\nu }\varepsilon
  ^{-2}\Ex_\vartheta \left|m_\varepsilon^\star \left(t\right)-m\left(\vartheta
  ,t\right)\right|^2=\sup_{\left|\vartheta -\vartheta _0\right|<\nu
  }\frac{\dot y\left(\vartheta ,t\right)^2 }{{\rm I}^t\left(\vartheta \right)}
  \left(1+o\left(1\right)\right),\\
&\sup_{\left|\vartheta -\vartheta
    _0\right|<\nu }\frac{\dot y\left(\vartheta ,t\right)^2 }{{\rm
      I}^t\left(\vartheta \right)}\longrightarrow \frac{\dot
    y\left(\vartheta_0 ,t\right)^2 }{{\rm
      I}^t\left(\vartheta_0 \right)},
\end{align*}
where the last limit corresponds to $\nu \rightarrow 0$. Of course this is not a proof.
The estimator $m_\varepsilon^\star \left(t\right)$  defined by  \eqref{29} is not exactly $m\left(\vartheta _{t,\varepsilon
}^\star, t\right) $ because if we substitute $\vartheta _{t,\varepsilon
}^\star $ in $m\left(\vartheta ,t\right)$, then we have to change the equality \eqref{27}.

Let us consider the difference between \eqref{29} and \eqref{30}. Below we use
Taylor expansions.
\begin{align*}
&m_\varepsilon^\star \left(t\right)-m\left(\vartheta ,t\right) =y_0\left[
N\left(\vartheta _{t,\varepsilon }^\star,t\right)-N\left(\vartheta ,t\right)\right]\\
&\qquad\qquad \qquad +\left[N\left(\vartheta
_{t,\varepsilon }^\star,t\right) -N\left(\vartheta,t\right)  \right] G\left(\vartheta _{t,\varepsilon }^\star
,X,t\right)\\
&\qquad\qquad \qquad  +N\left(\vartheta,t\right)  \left[ G\left(\vartheta _{t,\varepsilon }^\star
,X,t\right)- G\left(\vartheta ,X,t\right)\right]\\
&\qquad  =\left(\vartheta _{t,\varepsilon }^\star-\vartheta \right)
\left[y_0+G\left(\vartheta _{t,\varepsilon }^\star
,X,t\right) \right]\dot N(\tilde\vartheta _{t,\varepsilon }^\star,t)\\
&\qquad \qquad \qquad +\left(\vartheta _{t,\varepsilon }^\star-\vartheta \right)
N\left(\vartheta,t\right)  \dot G(\tilde\vartheta _{t,\varepsilon }^\star
,X,t)\\
&\qquad  =\left(\vartheta _{t,\varepsilon }^\star-\vartheta \right)
\left[\left(y_0+G\left(\vartheta
,X,t\right)\right) \dot N(\vartheta,t)+
N\left(\vartheta,t\right)  \dot G(\vartheta
,X,t)\right]\left(1+r_\varepsilon \right)\\
&\qquad  =\left(\vartheta _{t,\varepsilon }^\star-\vartheta \right)
\left[\left(y_0+G\left(\vartheta
,x,t\right) \right)\dot N(\vartheta,t)+
N\left(\vartheta,t\right)  \dot G(\vartheta
,x,t)\right]\left(1+\tilde R_\varepsilon \right)\\
&\qquad  =\left(\vartheta _{t,\varepsilon }^\star-\vartheta \right)
\dot y\left(\vartheta ,t\right)\left(1+\tilde R_\varepsilon \right).
\end{align*}
The conditions $C$ allow us to verify that the derivatives $\dot
N(\vartheta,t),\ddot N(\vartheta,t) $
are  bounded.  The derivatives of $G\left(\vartheta ,X,t\right)$  are
\begin{align*}
\dot G\left(\vartheta ,X,t\right)&=X_t\dot F\left(\vartheta ,t\right)-X_{\tau
  _\varepsilon }\dot F\left(\vartheta ,\tau
  _\varepsilon \right)-\int_{\tau _\varepsilon }^{t} \dot F'\left(\vartheta
,s\right)X_s{\rm d}s ,\\
\ddot G\left(\vartheta ,X,t\right)&=X_t\ddot F\left(\vartheta ,t\right)-X_{\tau
  _\varepsilon }\ddot F\left(\vartheta ,\tau
  _\varepsilon \right)-\int_{\tau _\varepsilon }^{t} \ddot F'\left(\vartheta
,s\right)X_s{\rm d}s .
\end{align*}
Therefore we have  for any $p>0$
\begin{align*}
\sup_{\left|\vartheta -\vartheta _0\right|<\nu }\Ex_\vartheta \left|\tilde
R_\varepsilon \right|^p\rightarrow 0
\end{align*}
as $\varepsilon \rightarrow 0$.

\section{Two Examples}

We consider two examples of construction of preliminary estimators.

{\bf Example 1.} In this example we have  the following system:
\begin{align*}
&{\rm d}X_t=\vartheta f_tY_t\,{\rm d}t+\varepsilon\sigma_t\,{\rm d}W_t,\quad
  X_0=0,\quad  0\leq t\leq T,\\
&{\rm d}Y_t=a_tY_t\,{\rm d}t+\varepsilon b_t\,{\rm d}V_t,\quad y_0\neq 0.
\end{align*}
For simplicity of exposition we suppose that functions $f_t,a_t,\sigma _t,b_t
, t\in \left[0,T\right]$ are positive
and bounded. The function $f_t\in {\cal C}_t^{\left(1\right)}$.  The
observations are $X^T=\left(X_t,0\leq t\leq T\right)$ and the
process $Y_t,0\leq t\leq T$ is {\it hidden}. The unknown parameter $\vartheta
\in\Theta =\left(\alpha ,\beta \right)$, $\alpha >0$.

The limit ($\varepsilon =0$) system is
\begin{align*}
x_\tau\left(\vartheta\right) &=\vartheta\int_{0}^{\tau }f_sy_s\,{\rm d}s .
\end{align*}
We put $\tau _\varepsilon =\varepsilon ^\delta $ and define   preliminary estimator  by the relation
\begin{align*}
\bar\vartheta _{\tau _\varepsilon }&=\frac{X_{\tau _\varepsilon
}}{\int_{0}^{\tau_\varepsilon  }f_sy_s\,{\rm d}s}=\frac{x_{\tau _\varepsilon
}\left(\vartheta _0\right)+\varepsilon \eta _{\tau _\varepsilon
}}{\int_{0}^{\tau_\varepsilon  }f_sy_s\,{\rm d}s}=\vartheta _0+ \frac{\varepsilon\eta _{\tau _\varepsilon
} }{\int_{0}^{\tau_\varepsilon  }f_sy_s\,{\rm d}s} .
\end{align*}
Hence
\begin{align*}
\Ex_{\vartheta _0}\left(\bar\vartheta _{\tau _\varepsilon
}-\vartheta_0\right)^2 =\varepsilon ^2\left(\int_{0}^{\tau_\varepsilon }f_sy_s\,{\rm d}s
\right)^{-2}\Ex_{\vartheta _0}\eta _{\tau _\varepsilon}^2\leq C\frac{\varepsilon ^2}{\tau _\varepsilon }.
\end{align*}

{\bf Example 2.} In the second example the partially observed system is
\begin{align*}
&{\rm d}X_t= f_tY_t\,{\rm d}t+\varepsilon\sigma_t\,{\rm d}W_t,\quad
  X_0=0,\quad  0\leq t\leq T,\\
&{\rm d}Y_t=\vartheta a_tY_t\,{\rm d}t+\varepsilon b_t\,{\rm d}V_t,\quad
  y_0\neq 0.
\end{align*}
We have $Y_t=y_t\left(\vartheta _0\right)+\varepsilon \xi _t$ and   the limit equation is
\begin{align*}
y_t\left(\vartheta _0\right)=y_0+\vartheta_0 \int_{0}^{t}a_sy_s\left(\vartheta _0\right)\,{\rm d}s .
\end{align*}
Therefore  we can write
\begin{align*}
X_t&=\int_{0}^{t}f_sY_s\,{\rm d}s+\varepsilon \int_{0}^{t}\sigma _s\,{\rm
  d}W_s\\ &=\int_{0}^{t}f_sy_s\left(\vartheta _0\right)\,{\rm d}s+\varepsilon
\int_{0}^{t}f_s\xi _s\,{\rm d}s+\varepsilon \int_{0}^{t}\sigma _s\,{\rm
  d}W_s\\
&=y_0\int_{0}^{t}f_s{\rm d}s+\vartheta
_0\int_{0}^{t}f_s\int_{0}^{s}a_ry_r\left(\vartheta _0\right){\rm d}r{\rm d}s +\varepsilon
\int_{0}^{t}f_s\xi _s{\rm d}s+\varepsilon \int_{0}^{t}\sigma _s\,{\rm
  d}W_s.
\end{align*}
For the small values of $\tau $ we have
\begin{align*}
X_\tau-y_0\int_{0}^{\tau}f_s\,{\rm d}s&=\vartheta _0\, f_0\,a_0\,y_0\,
\frac{\tau ^2}{2}+O\left(\tau ^3\right)+\varepsilon
\int_{0}^{\tau}f_s\xi _s{\rm d}s+\varepsilon \int_{0}^{\tau}\sigma _s\,{\rm
  d}W_s  .
\end{align*}
The preliminary estimator can be defined as follows
\begin{align*}
\bar\vartheta _{\tau _\varepsilon }= \frac{2}{f_0\,a_0\,y_0\,\tau _\varepsilon
  ^2}\left(X_{\tau _\varepsilon }-y_0\int_{0}^{\tau_\varepsilon }f_s\,{\rm
  d}s\right).
\end{align*}
The error of estimation is
\begin{align*}
\Ex_{\vartheta _0}\left(\bar\vartheta _{\tau _\varepsilon }-\vartheta _0
\right)^2\leq C \tau _\varepsilon ^2+C\frac{\varepsilon ^2 }{\tau _\varepsilon
  ^3}=C\left(\varepsilon ^{2\delta } +\varepsilon ^{2-3\delta }\right).
\end{align*}
Hence if $\delta \in \left(0, 2/3\right)$, then this estimator is
consistent. The optimal choice, which minimizes the error is $\delta =2/5$.

\section{Discussion}

The main conditions in the construction of One-step MLE $\vartheta _\varepsilon ^\star$ and One-step
MLE-process $\vartheta _{t,\varepsilon} ^\star, \tau _\varepsilon <t\leq T$  are $\left|\dot f\left(\vartheta ,0\right)\right|>0$ or $\left|\dot
a\left(\vartheta ,0\right)\right|>0$. Consider the situation where the both
conditions are not fulfilled. Suppose that there exists $\tau _0\in
\left(0,T\right)$ such that
$$
f\left(\vartheta ,t\right)=f\left(t\right)\1_{\left\{0\leq t<\tau _0\right\}}
+\vartheta g\left(t\right)\1_{\left\{\tau _0\leq  t \leq T\right\}},
$$
 and $a\left(\vartheta ,t\right)=a\left(t\right), 0\leq
t\leq T $. Then the construction of these estimators can be done as
follows. First we solve the filtration equations on the interval $\left[0,\tau
  _0\right]$
\begin{align*}
{\rm d}m\left(t\right)&=a\left(t\right)m\left(t\right){\rm d}t+\frac{\gamma_*
  \left(t\right)f\left(t\right) }{\sigma \left(t\right)^2} \left[{\rm
    d}X_t-f\left(t\right)m\left(t\right){\rm d}t \right],\quad
m\left(0\right)=y_0,\\ \frac{\partial \gamma_* \left(t\right)}{\partial
  t}&=2a\left(t\right)\gamma_* \left(t\right)-\frac{\gamma_*
  \left(t\right)^2f\left(t\right)^2}{\sigma
  \left(t\right)^2}+b\left(t\right)^2,\quad \gamma_* \left(0\right)=0.
\end{align*}
Then we can consider $m\left(\tau _0\right),\gamma_* \left(\tau _0\right)$ as
initial values for the filtration equations on the interval $\left[\tau
  _0,T\right]$ and the preliminary estimator $\bar\vartheta _{\tau _0+\tau
  _\varepsilon }$ can be constructed by observations $X_{\tau _0}^{\tau
  _\varepsilon }=\left(X_t, \tau _0\leq t\leq \tau
  _0+\tau _\varepsilon \right)$.

\bigskip

Another question concerns the length of the ``learning interval'' $\left[0,\tau
  _\varepsilon \right]$, where  $\tau _\varepsilon =\varepsilon ^\delta $ and
$\delta \in \left(0,1\right)$.  If we need the approximation of
$m\left(\vartheta _0,t\right)$ for the values $t$ smaller than the given  $\tau
_\varepsilon $, then we can use the construction of  Two-step MLE-process
described in the work  \cite{Kut17}.

\bigskip

In the case of multidimensional $\vartheta \in \Theta \subset {\cal R}^d$ we
have to suppose that the observed process is as well multidimensional
$X^T=\left(X_1^T, \ldots, X_k^T\right)$, where $X_j^T=\left(X_{j,t}, 0\leq
t\leq T\right)$ and $k\geq d$. For example,
\begin{align*}
{\rm d}X_{j,t}=f_j\left(\vartheta ,t\right)Y_t\,{\rm d}t+\varepsilon \sigma
_j\left(t\right)\,{\rm d} W_j\left(s\right),\quad X_{j,0}=0,\quad j=1,\ldots,k.
\end{align*}

\bigskip

{\bf Acknowledgment.}  The authors gratefully acknowledge financial support from the National Natural Science Foundation of China (NSFC) under grant 11701329.


\begin{thebibliography}{99}

\bibitem{A83} Arato, M. (1983) {\it  Linear Stochastic Systems with Constant
  Coefficients A Statistical Approach.} Lecture Notes in Control and
  Inform. Sci., 45, New York: Springer-Verlag.

\bibitem {BRR98} Bickel, P.J., Ritov, Y. and Ryd\'en, T. (1998) Asymptotic
  normality of the maximum likelihood estimator for general hidden Markov
  models. {\it   Ann. Statist.}, 26, 4, 1614-1635.

\bibitem {CMT05} Capp\'e, O., Moulines, E. and Ryd\'en, T. (2005) {\it  Inference
  in Hidden Markov Models}. Springer, N.Y.


\bibitem {PCh09} Chigansky, P. (2009) Maximum likelihood estimation for hidden
  Markov models in continuous time. {\it   Statist. Inference Stoch. Processes},
  12, 2, 139-163.


\bibitem {EAM95} Elliott, R.J., Aggoun, L. and Moor, J.B. (1995) \textsl{Hidden
  Markov Models.} Springer,  N.Y.

\bibitem {EM02} Ephraim, Y., Mehrav, N. (2002) Hidden Markov
  processes. {\it  IEEE Trans. Inform. Theory,} 48, 6, 1518-1569.

\bibitem{F25} Fisher, R.A. (1925) {\it Theory of Statistical Estimation.}
  Proc. Camb. Phil. Soc., 22, 700-725.

\bibitem{GL95} Gill, R.D. and Levit, B.Ya  (1995) Application of the van Trees
  inequality: a Bayesian Cramer-Rao bound. {\it Bernoulli}, 1, 59-79.

\bibitem{G84} Golubev, G.K. (1984) Fisher's method of scoring in the problem of
  frequency estimation. {\it  J. of Soviet Math.,} 25, 3, 1125-1139.

\bibitem{Gus00} Gustafsson, F. (2000) {\it Adaptive Filtering and Change
  Detection.} J. Wiley\&Sons, N.Y.

\bibitem{HMS03} Hide, C., Moore, T. and Smith, M. (2003) Adaptive Kalman
  filtering for low cost ING/GPS. {\it The Journal of Navigation}, 56, 1, 143-152.


\bibitem{HCCL03} Hu, C., Chen, W., Chen, Y. and Liu, D. (2003) Adaptive Kalman
  filtering for vehicule navigation. {\it J. Global Positioning systems}, 2,
  1, 42-47.

\bibitem{K-B61} Kalman, R.E, and Bucy, R.,S. (1961) New results in linear
  filtering and prediction theory. {\it  Trans. ASME}, 83D, 95-100.

\bibitem{KU15} Kamatani, K. and Uchida, M. (2015) Hybrid multi-step
  estimators for stochastic differential equations based on sampled data.
 {\it  Statist. Inference  Stoch. Processes.} 18, 2, 177-204.

\bibitem{Kh05} Khasminskii, R.  (2005) Nonlinear filtering of smooth
  signals. {\sl Stochastics and Dynamics,} 5, 1, 27-35.



\bibitem{KhK18} Khasminskii, R. Z. and Kutoyants, Yu. A.  (2018) On parameter
  estimation of hidden telegraph process.
  {\it Bernoulli}, 24, 3, 2064-2090.
	
\bibitem{Kut84} Kutoyants, Y.A. (1984) {\it  Parameter Estimation for
  Stochastic Processes.} Heldermann, Berlin.

\bibitem{Kut94} Kutoyants, Y.A. (1994) {\it  Identification of Dynamical
  Systems with Small Noise.} Kluwer Academic Publisher, Dordrecht.

\bibitem{Kut04} Kutoyants, Yu.A. (2004) {\it  Statistical Inference for
  Ergodic Diffusion Processes.} Springer, London.

\bibitem{Kut17} Kutoyants, Yu.A. (2017) On the multi-step MLE-process for
  ergodic diffusion.  {\it  Stochastic     Process. Appl.}, 127, 2243-2261.

\bibitem{K19} Kutoyants, Yu. A.  (2019)  On parameter estimation of the hidden
  Ornstein--Uhlenbeck process. {\it   J. Multivar. Analysis}, 169, 1,
  248-269.

\bibitem{Kut19} Kutoyants, Yu. A.  (2019)  On parameter estimation of the
  hidden ergodic   Ornstein--Uhlenbeck process. Submitted. (arXiv:1902.08500)

\bibitem{KM16} Kutoyants, Yu. A. and Motrunich, A.  (2016) On multi-step
  MLE-process for Markov sequences. {\it   Metrika}, 79, 705-724.

\bibitem {KZ14} Kutoyants, Y.A. and Zhou, L. (2014) On approximation of the
  backward stochastic differential equation.
  {\it J. Stat. Plann. Infer. } 150, 111-123.

\bibitem{LS} Liptser, R.S. and  Shiryayev, A.N.  (2001) {\it  Statistics of
  Random Processes, I. General Theory.} 2nd Ed., Springer, N.Y.

\bibitem{Ru91} Rutan, S.C. (1991) Adaptive Kalman filtering. {\it
  Anal. Chem.}, 63 (22), 1103A-1109A.


\end{thebibliography}
\end{document}